\documentclass[journal]{IEEEtran}
\usepackage{amsfonts}
\usepackage{bm}
\usepackage{amsmath,amsfonts,amsthm,amssymb}
\usepackage{setspace}
\usepackage{Tabbing}
\usepackage{fancyhdr}
\usepackage{lastpage}
\usepackage{extramarks}
\usepackage{chngpage}
\usepackage{algorithmic}
\usepackage{soul,color}
\usepackage{epsfig}
\usepackage{graphicx,float,wrapfig,subfigure}
\usepackage{dsfont}
\usepackage{longtable}
\usepackage{wrapfig}
\usepackage{stfloats}
\usepackage{cite}
\usepackage{graphicx}
\usepackage{array}
\usepackage{multirow}{\tiny }
\usepackage{multicol}
\usepackage{colortbl}
\usepackage{tabularx}
\usepackage{mdwmath}
\usepackage{mdwtab}
\usepackage{color}
\usepackage{verbatim}
\usepackage{amsmath}
\usepackage{tikz}
\usetikzlibrary{calc}
\usepackage{flowchart}
\usetikzlibrary{shapes.geometric, arrows}
\usepackage{overpic}
\usepackage{epstopdf}
\usepackage{url}
\usepackage[colorlinks,linkcolor=black,anchorcolor=black,citecolor=black,urlcolor=black]{hyperref} 
\usepackage[flushleft]{threeparttable}
\allowdisplaybreaks

\begin{document}

\title{
	A Second Order Cone Programming Model for Planning PEV Fast-Charging Stations
	}

\author{
Hongcai~Zhang,~\IEEEmembership{Student Member,~IEEE,}
Scott~J.~Moura,~\IEEEmembership{Member,~IEEE,}
Zechun~Hu,~\IEEEmembership{Senior Member,~IEEE,}
Wei~Qi,
and~Yonghua~Song,~\IEEEmembership{Fellow,~IEEE}

\thanks{
	This work was supported in part by the National Key Research and Development Program (2016YFB0900103) and the National Natural Science Foundation of China (51477082).
	
	H. Zhang, Z. Hu and Y. Song are with the Department of Electrical Engineering, Tsinghua University, Beijing, 100084, P.~R.~China (email:  zhang-hc13@mails.tsinghua.edu.cn; zechhu@tsinghua.edu.cn; yhsong@tsinghua.edu.cn).
	
	S. J. Moura is with the Department of Civil and Environmental Engineering, University of California, Berkeley, California, 94720, USA, and also with the Smart Grid and Renewable Energy Laboratory, Tsinghua-Berkeley Shenzhen Institute, Shenzhen, 518055, P.~R.~China (email:  smoura@berkeley.edu).

	W. Qi is with the Desautels Faculty of Management, McGill University, Montreal, Quebec, Canada H3A 1G5 (email:  qiwei.0216@gmail.com). 
	}
}

\maketitle

\begin{abstract} 
	This paper studies siting and sizing of plug-in electric vehicle (PEV) fast-charging stations on coupled transportation and power networks. 
	We develop a closed-form model for PEV fast-charging stations' service abilities, which considers heterogeneous PEV driving ranges and charging demands.
	We utilize a modified capacitated flow refueling location model based on sub-paths (CFRLM\_SP) to explicitly capture time-varying PEV charging demands on the transportation network under driving range constraints. We explore extra constraints of the CFRLM\_SP to enhance model accuracy and computational efficiency.
	We then propose a stochastic mixed-integer second order cone programming (SOCP) model for PEV fast-charging station planning. The model considers the transportation network constraints of CFRLM\_SP and the power network constraints with AC power flow. 
	Numerical experiments are conducted to illustrate the effectiveness of the proposed method.
\end{abstract}


\section*{Nomenclature}
\begin{table}[H]
	\renewcommand{\arraystretch}{1.05}
	\begin{normalsize}
		\begin{tabular}{ p{.06\textwidth}p{.39\textwidth}}
			\multicolumn{2}{l}{\textit{Indices/Sets}}\\
			$e_{(k)}$& Index of (type $k$) plug-in electric vehicles (PEVs).\\
			$i/\mathcal{I}_{(o/q)}$&Index/set of transportation nodes (on sub-path $o$/path $q$), $i\in \mathcal{I}_{(o/q)}$.\\
			$\mathcal{I}_{m}$& Set of transportation nodes connected to distribution bus $m$.\\
			$k/\mathcal{K}$& Index/set of PEV types, $k \in \mathcal{K}$.\\
			$m/n/h$&Index of buses of the distribution network. $m/n/h \in\mathcal{M}$. For the reference bus, $m/n/h=0$.\\
			$(m,n)/\mathcal{L}$&Index/set of lines of the distribution network. $(m,n)$ is in the order of bus $m$ to bus $n$, i.e., $m\rightarrow n$, and bus $n$ lies between bus $m$ and bus 0. $(m,n) \in \mathcal{L}$.\\
		\end{tabular}
	\end{normalsize}
\end{table}
\begin{table}[H]
	\renewcommand{\arraystretch}{1.05}
	\begin{normalsize}
		\begin{tabular}{ p{.06\textwidth}p{.39\textwidth}}
		$\mathcal{M}_{(\rightarrow m)}$&Set of buses of the distribution network (that are connected to bus $m$ and bus $m$ lies between them and bus 0).\\
			$o/\mathcal{O}_{(q,k)}$&Index/set of sub-paths (of PEV type $k$ on path $q$), $o \in \mathcal{O}_{(q,k)}$.\\
			$q/\mathcal{Q}_{(i)}$&Index/set of paths (travel by node $i$), $q \in \mathcal{Q}_{(i)}$. Each path corresponds to one OD pair.\\
			$t$ & Index of time intervals.\\
			$\omega/\Omega$& Index/set of scenarios. $\omega \in \Omega$.\\
			\multicolumn{2}{l}{\textit{Functions}}\\
			$f(\cdot)$&Probability density function (PDF) of a normal distribution.\\
			$F(\cdot)$&Cumulative distribution function (CDF) of a normal distribution.\\
			$\Phi(\cdot)$&CDF of the standard normal distribution.\\
			\multicolumn{2}{l}{\textit{Parameters of the PEVs}}\\
			$R_k$& Driving range of a type $k$ PEV, in km.\\
			$t_{e}^\text{a/c}$&Arrival (charge start) time of PEV $e$ in a station. $t_{0}^\text{a}=0$.\\
			$t_{e}^\text{d}$&Departure (charge end) time of PEV $e$ in a station.\\
			$T_{(k)}$& Required charging time for a (type $k$) PEV, in h.\\
			$y_{(k)}^\text{ev}$& The arrival number of (type $k$) PEVs in a charging station.\\
			\multicolumn{2}{l}{\textit{Parameters of the transportation network}}\\
			$d_{i,j}$ & Distance between node $i$ and $j$, in km.\\
			$\lambda_{q(,k)}$&Volume of (type $k$) PEV traffic flow on path $q$, in $\text{h}^{-1}$.\\
			$\lambda_{q,i,k,\omega,t}$&Volume of type $k$ PEV traffic flow (Poisson parameter) on path $q$, at node $i$, during time $t$, in scenario $\omega$, in $\text{h}^{-1}$.\\
			\multicolumn{2}{l}{\textit{Parameters of the planning model}}\\
			$\alpha$& Service level of charging stations.\\
			$c_{1,i}$&Fixed costs for building a new charging station at node $i$, including buildings costs etc, in \$.\\
			$c_{2,i}$&Costs for adding an extra charging spot in a station at node $i$, including land use costs, spot purchase costs etc., in \$.\\
			$c_{3,i}$&Per-unit cost for distribution line at $i$, in \$/(kVA$\cdot$km).\\
			$c_{4,i}$&Per-unit cost for substation capacity expansion at $i$, in \$/kVA.\\
			$c_\text{e}$&Per-unit cost for energy purchase, in \$/kWh.\\
			$c_\text{p}$&Per-unit penalty costs for unsatisfied PEV power, in \$/kWh.\\
			$\Delta t$&Time interval, representing one hour in this paper.\\
		\end{tabular}
	\end{normalsize}
\end{table}

\begin{table}[H]
	\renewcommand{\arraystretch}{1.05}
	\begin{normalsize}
		\begin{tabular}{ p{.06\textwidth}p{.39\textwidth}}
		$\overline{I_{mn}}$&Upper limit of branch current of line $(m,n)$, in kA.\\
			$l_{i}$&Required distribution line length to install a charging station at transportation node $i$, in km.\\
			$\zeta$& Capital recovery factor, which converts the present investment costs into a stream of equal annual payments over the specified time of $Y$ at the given discount rate $r$. $\zeta=\left({r(1+r)^{Y}}\right) / \left({(1+r)_{}^{Y}-1}\right)$. $Y$ is the service live of the charging stations, in year.\\
			$\pi_\omega$& Probability of scenario $\omega$, in \%. \\
			$p^{\text{spot}} $&Rated charging power of a charging spot, in kW.\\
			$P_{i,0}^\text{sub}$&Initial substation capacity available at node $i$, in kVA.\\
			$s_{m,\omega,t}^\text{b}$&Apparent base load at bus $m$, in kVA.\\
			$\underline{V_{m}}/\overline{V_{m}}$&Lower/upper limit of nodal voltage at bus $m$, in kV.\\
			$\overline{y}_{i}$&Maximum number of charging spots located at node $i$.\\
			$z_{mn}$& Impedance of branch $(m,n)$, in ohm. $z_{mn}^*$ is its conjugate.\\
			\multicolumn{2}{l}{\textit{Decision Variables}}\\
			$\lambda_{(i,k,\omega,t)}$&Volume of (type $k$) PEVs that need charging (at node $i$, during $t$, in scenario $\omega$), in $\text{h}^{-1}$.\\
			$\gamma_{q,i(,k)}$&Charge choice of (type $k$) PEVs on path $q$ at node $i$: $\gamma_{q,i(,k)}=1$, if PEVs get charged; $\gamma_{q,i(,k)}=0$, otherwise.\\
			$l_{mn,\omega,t}$& Square of the magnitude of line $(m,n)$'s current during $t$ in scenario $\omega$, in $\text{kA}^2$.\\
			$p_{i,\omega,t}^\text{ev}$&Active PEV power at node $i$ during $t$ in scenario $\omega$, in kW.\\
			$p_{\text{un},(i),\omega,t}^\text{ev}$&Unsatisfied active PEV charging power (at node $i$) during $t$ in scenario $\omega$, in kW.\\
			$p_{m,\omega,t}$&Total active power injection at bus $m$ during $t$ in scenario $\omega$, in kW.\\
			$P_i^\text{sub}$&Substation capacity expansion at node $i$, in kVA.\\
			$s_{m,\omega,t}$&Total apparent power injection at bus $m$ during $t$ in scenario $\omega$, in kVA. For a distribution system, $s_{0,\omega,t}$ (at bus $0$) is also the power consumption of the whole distribution system.\\
			$s_{m,\omega,t}^\text{ev}$& Apparent PEV power at bus $m$ during $t$ in scenario $\omega$, in kVA.\\
			$S_{mn,\omega,t}$&Apparent power flow in line $(m,n)$ (from bus $m$ to bus $n$) during $t$ in scenario $\omega$, in kVA.\\
			$v_{m,\omega,t}$&Square of the nodal voltage at bus $m$ during $t$ in scenario $\omega$, in $\text{kV}^2$.\\
			$x_i$&Charging station investment decision at node $i$: $x_i=1$, if there is a station at node $i$; $x_i=0$, otherwise.\\
			$y_i^\text{cs}$&Number of invested charging spots (at node $i$).\\
			$z_{o_k,i}$&Charge choice of PEVs on sub-path $o_k$ at node $i$: $z_{o_k,i}=1$, if charged; $z_{o_k,i}=0$, otherwise.\\
		\end{tabular}
	\end{normalsize}
\end{table}

\section{Introduction}
\IEEEPARstart{P}{lug-in} electric vehicles (PEVs) are regarded as a promising tool to promote energy sustainability and combat climate change. \textcolor{black}{Compared with traditional  internal combustion engine vehicles, PEVs have higher energy efficiency. 
According to a report from the Argonne National Laboratory \cite{EV_Emission_Argonne2010}, the electrification of transportation can significantly reduce petroleum energy use and help to relief the pressure of energy crisis. The emissions of PEVs strongly depend on their electricity generation mix for recharging. With the fast development of sustainable resources, e.g., wind, photovoltaic and hydro power, adopting PEVs can also significantly reduce the global greenhouse gas emissions\cite{EV_Emission_Argonne2010}. Furthermore, because PEV charging demands are usually flexible, they may provide various power grid services under proper management, e.g., load valley filling\cite{ma2013decentralized,liu2016aggregation,li2017data}, participating energy market\cite{Sarker2016,neyestani2017plug}, promoting renewable power adoption \cite{Wang2015,Vaya2016}, or providing ancillary services \cite{Zhang_V2GTPS2016,Wu2016,Liuhui_EVregulation_2016}. This will also enhance the sustainability and low-carbonization of future power systems. }

\textcolor{black}{Because of the aforementioned advantages of PEVs, governments around the world have devoted great efforts to promoting their development. 
As a result, the PEV market has experienced an explosive growth in recent years. Over 774,000 PEVs were sold worldwide in 2016, which increased by 42\% compared to 2015\cite{EV_Population}. A report by Bloomberg forecasted that the sales of PEVs will hit 41 million by 2040 worldwide, accounting for 35\% of new light duty vehicle sales\cite{EV_Forecast}.}

Well developed infrastructure for PEV charging is the prerequisite for the promotion of PEV adoption.
The growing PEV population is leading to massive investments in charging infrastructure. For example, in China, 4.8 million distributed charging spots and more than twelve thousand fast-charging stations are planned for construction by 2020\cite{EVStationPlan_chinadaily}. PEV charging facilities generally fall into two categories: 1) distributed charging spots with slow (or normal) power chargers; 2) centralized fast-charging stations with high power chargers\cite{Plan_Zhang2015}. In urban areas, the distributed charging spots, located in private or public parking lots etc., serve as the primary means for PEV charging; while on intercity corridors, the fast-charging stations are the major charging infrastructure for PEVs.
Driven by the urgent demands from industry, the planning of both kinds of PEV charging facilities have become important research focuses \cite{Plan_T_Dong2016,Plan_Zhang2016,Plan_Zhang2015}.
{This paper studies siting and sizing of fast-charging stations on coupled transportation and power networks, solving three major sub-problems:
\begin{enumerate}
	\item How to model a fast-charging station's service ability considering heterogeneous PEV charging demands? In other words, when heterogeneous PEV charging demands are given at a charging station, how many spots should be installed there to offer adequate charging service quality?
	\item How to model PEVs' charging behaviors on transportation networks? In other words, where should PEVs get charged (or charging stations be located) so that their driving demands can be satisfied? 
	\item PEV charging stations are the coupling points of the transportation and power networks. How to describe the two networks' coupling relationship so that we can optimize the sites and sizes of fast-charging stations subjected to their coupled constraints? 
\end{enumerate}}

{Modeling one single fast-charging station's service ability, i.e., the number of demands that a station with a certain number of charging spots can satisfy, is the first step for the planning. PEVs' long and heterogeneous charging time makes the modeling of charging stations different from that of traditional gasoline stations.} Some previous literature\cite{FRLM_Kuby2005,FRLM_MirHassani2013,FRLM_Chung2015} studied un-capacitated charging station planning, which assumed infinite service ability for each station. Some others (e.g., \cite{CFRLM_Upchurch2009}) approximated a station's service ability to be proportional to its facility number. This linear model ignores the ``scale effect" of stations' service abilities, i.e., when the facility number in a station grows, the average service ability of a single facility increases because of the randomness of demands. To handle the ``scale effect," queuing theory was widely adopted, e.g., in \cite{Queuing_Li2012,Plan_TE_Yao2014,Plan_TE_Xiang2016,Queuing_Fan2015}. \textcolor{black}{Though queuing models provide a comparatively precise way to model a PEV charging station's service ability, they have no closed-form formulations. Users have to solve difficult nonlinear optimization problems to obtain the required number of charging spots for a certain number of demands, which makes it hard to be directly implementable in the planning model. Most of the published papers that have adopted queuing models applied heuristic algorithms to solve their problems\cite{Queuing_Li2012,Plan_TE_Yao2014,Plan_TE_Xiang2016,Queuing_Fan2015}. Besides, because the battery capacities of PEVs on the market are heterogeneous, their required service durations in a charging station may have heterogeneous distributions, which is not considered in the aforementioned models.}

{The limited driving range is another key characteristic that distinguishes PEVs from traditional internal combustion engine vehicles. Unlike gasoline stations, planning PEV fast-charging stations should consider PEVs' driving range constraints, which determine where and when their charging demands may happen on transportation networks.} Two different ways to explicitly model PEVs' driving range constraints proposed in published literature are respectively named in this paper as: 1) the flow refueling location model based on expanded networks (FRLM\_EN); 2) the flow refueling location model based on sub-paths (FRLM\_SP). Both models use origin-destination (OD) traffic flow to estimate charging demands. FRLM\_EN was first proposed in \cite{FRLM_Kuby2005} and then reformulated and simplified in \cite{FRLM_MirHassani2013,FRLM_Chung2015}. FRLM\_SP was developed in \cite{Plan_T_Mak2013} for battery swapping station planning. In both FRLM\_EN and FRLM\_SP, only  peak-hour OD traffic flow is considered.  However, modeling time-varying traffic flow is necessary to evaluate the impact of PEV charging on dynamic power networks.

{PEV charging stations are the coupling points of transportation and power networks. Their investments and operations are constrained by both networks. Hence, planning PEV charging stations should consider the coupling relationship between the transportation and the power networks.} Very few published papers study this coupling. 
In \cite{Plan_TE_Wang2013, Plan_TE_Yao2014, Plan_TE_Xiang2016, Plan_TE_Sadeghi2014, Plan_TE_Luo2015,Plan_TE_He2013}, detailed power network constraints were considered. However, references \cite{Plan_TE_Wang2013, Plan_TE_Yao2014,Plan_TE_Xiang2016} modeled the transportation network without explicitly considering driving range constraints and they considered low voltage distribution networks with service radiuses much smaller than a typical PEV's driving range. Reference \cite{Plan_TE_Sadeghi2014} assumed the PEV charging demands to be uniformly distributed across the target area. In \cite{Plan_TE_Luo2015}, the authors studied charging station siting problem based on game-theoretical modeling and simulation. The impact of PEV charging on the power grid was assumed to be proportional to the charging power. In \cite{Plan_TE_He2013}, an equilibrium modeling framework for PEV charging station planning in a coupled transportation and power network is proposed. The authors assumed the transmission nodal electricity prices will influence PEV charging choices and the traffic flow. In practice, the geographical distance between two transmission nodes and the costs for a PEV to travel from one node to another is high so that the nodal prices effect may be insignificant.


In our previous work \cite{Zhang_PlanFRLMTSG2016}, we proposed a mixed integer linear programming model for PEV fast-charging station planning on coupled transportation and power networks. We used queuing theory to model charging stations' service abilities for PEVs with homogeneous driving ranges and adopted piecewise linearization to retain linearity (at the price of introducing additional binary variables). We utilized capacitated FRLM\_EN (CFRLM\_EN) to model the transportation network and only considered peak-hour OD traffic flow. Kirchhoff’s Law was utilized to roughly approximate the electrical constraints of distribution networks.

Compared with the published literature including our previous work \cite{Zhang_PlanFRLMTSG2016}, the major contributions of this paper include:
\begin{enumerate}
	\item \textcolor{black}{We propose a service performance metric for PEV fast-charging stations, i.e., the \textit{service level}. It measures the probability that the charging demands arriving in a given time interval can be directly fulfilled without extra waiting time (caused by limited charging capacity). Then, we develop a closed-form model based on the service level metric to describe the service abilities of PEV fast-charging stations.} Compared with linear models and queuing models in published literature, the advantages of the proposed model are twofold:
	\begin{enumerate}
		\item It is in a simple closed form and can be modeled as an SOCP constraint so that it can be easily implemented in the planning model and solved by an off-the-shelf solver, which ensures optimality;
		\item It considers heterogeneous PEV driving ranges so that the modeling is more practical and accurate.
	\end{enumerate}
	\item A modified capacitated FRLM\_SP (CFRLM\_SP) is designed to explicitly capture time-varying PEV charging demands generated from dynamic OD traffic flow on the transportation network under driving range constraints. The advantages of this model are threefold:
	\begin{enumerate}
		\item Compared with the CFRLM\_EN, it introduces fewer binary decision variables so that the problem scale is reduced;
		\item We develop extra constraints based on practical operation analysis to further enhance  its computational efficiency;
		\item \textcolor{black}{Modeling dynamic OD traffic flow allows us to evaluate time-varying PEV charging power, which is crucial for the secure operation of distribution systems considering that the base load profiles are also time-varying.}
	\end{enumerate}
	\item A stochastic mixed-integer SOCP is developed for PEV charging station planning considering both the transportation network constraints and the power network constraints with AC power flow. By modeling AC power flow instead of its linear approximation in \cite{Zhang_PlanFRLMTSG2016}, we can evaluate the planning results' influence on important parameters of distribution system operations including active power losses and voltage deviations. Besides, we adopt convex SOCP relaxation of AC power flow so that its global optimal solution can be obtained\cite{OPF_tree_Exactness_Gan2015}. The planning model can be efficiently solved by the branch-and-cut method using an off-the-shelf solver.
\end{enumerate}
Besides, numerical experiments are conducted to validate the proposed planning method.

\textcolor{black}{Note that we study this problem from the perspective of a social planner with an objective to maximize the social welfare. We assume that the planner has access to both transportation and power system information. In scenarios when the power utility companies are also the charging station investors, e.g., State Grid Corporation of China, the proposed method is also applicable in practice. The targeted planning area in this paper is a highway transportation network powered by a high voltage distribution network with large service radius as in \cite{Zhang_PlanFRLMTSG2016}. The proposed method can be easily extended to scenarios when higher voltage level transmission systems are also covered in the targeted planning area.}

{The service performance metric and the service ability model are introduced in Section II. The latter is used to determine the sizes of fast-charging stations given the demands. Section III introduces the modified CFRLM\_SP, which determines the feasible set of the locations where PEVs get charged (or fast-charging stations are constructed). Section IV formulates the mixed-integer SOCP planning model which subjects to the aforementioned service ability model and the CFRLM\_SP. The AC power flow model and other constraints are also introduced. Case studies are described in Section V and Section VI concludes the paper.}

\section{Service Ability Model of Charging Stations}\label{sec:SRM}
{This section first proposes a service performance metric, i.e., the \textit{service level}, for PEV fast-charging stations. Then, based on the proposed metric, we develop a closed-form model to describe the service ability of a PEV fast-charging station, which is called the \textit{service level model} in this paper. This model is used to size a PEV fast-charging station, i.e., to determine the required number of charging spots in a station given the number of charging demands.} 

\subsection{The Service Level Metric}
Given the charging demands at a candidate location, the planer should determine the number of charging spots to be constructed there. To guarantee that the designed station can provide adequate service quality in the future, the planner should determine this number according to a proper charging service performance metric. 
This section introduces a novel performance metric, i.e., the \textit{service level}, for PEV fast-charging stations.
The “service level” is a popular performance metric in inventory management. It measures the probability (rate) that all customer orders arriving within a given time interval will be completely delivered from stock on hand, i.e. without delay\cite{SupplyChain_Max2011}. 
Reference \cite{Plan_T_Mak2013} used it to describe the service quality of a battery swapping station. 
Inspired by this metric in inventory management, we define the service level of a fast-charging station as follows:

{\textbf{Definition 1} \textit{(Service level):} The service level of a fast-charging station represents the probability that any PEV can be charged for at least its required units of time, i.e., $T_e$ for PEV $e$, without waiting in the station.

\textcolor{black}{Since the user experiences in a fast-charging station are strongly related to the users' waiting time before getting charged, the proposed service level metric can effectively describe the service quality of a fast-charging station.}

\textcolor{black}{Based on the service level metric, we propose that the planner determine the number of charging spots at a candidate location subject to the following service quality criterion:}

\textcolor{black}{\textbf{Criterion 1}: The service level of the fast-charging station is $\alpha$ or greater. Mathematically, $\text{Pr}( t_{e}^\text{c} = t_{e}^\text{a}~\&~t_{e}^\text{d}-t_{e}^\text{c}\geq T_e)\geq \alpha, ~ \forall e$.}

\textcolor{black}{By tuning the parameter $\alpha$ (which is called the service level criterion in this paper), the planner can effectively control the designed station's future service quality.} 

\textcolor{black}{When the charging demands are given, the minimum number of charging spots that should be installed in the station, i.e., $y^\text{cs}$, is a function of the service level criterion  $\alpha$. In the following subsections, we derive a closed-form approximation for this function (denoted by $y^\text{cs}(\alpha)$). We begin the derivations assuming that the PEVs have homogeneous driving ranges. Then, we extend the result for the scenarios when PEVs are heterogeneous. This function $y^\text{cs}(\alpha)$ defines the service ability of a fast-charging station subject to the service level criterion $\alpha$, i.e., the service level model. At the end of this section, we show how to reformulate this service level model into a mixed-integer SOCP.}

\subsection{The Service Level Model with Homogeneous PEVs}
\textcolor{black}{
We first consider the scenario when all the PEVs are homogeneous so that the required time to fully charge them with depleted batteries are the same. We denote the time by $T$, i.e., $T_e =T, ~\forall e$. 
To derive the closed-form approximation for the function $y^\text{cs}(\alpha)$, we introduce the following assumptions:}

\textbf{[A1]} The PEVs arrive in a fast-charging station following a Poisson process.

\textbf{[A2]} The PEVs are served based on a \textit{first-in first-out rule} and new arriving PEVs need not to wait, which means that, when all the charging spots are occupied and a new PEV arrives, the on-board PEV that has been charged the most must leave and spare its spot to the new one. Thus, a PEV $e$ may leave the station if it : 1) has got charged for $T_e$ units of time; 2) has to spare its spot to a new PEV. Therefore, with this assumption, we always have $t_{e}^\text{c} = t_{e}^\text{a}$.\footnote{This mild assumption is made for the convenience of modeling. In practice, a new arriving PEV may need to wait before one spot is spared for it.}

Let $y^\text{ev}$ denote the arrival number of PEVs in the station in a duration of $T$; $\lambda$ denote the Poisson arrival rate; then, based on assumption [A1], $y^\text{ev}$ follows a Poisson distribution with parameter $T\lambda$, i.e., $y^\text{ev}\sim Poisson(T\lambda)$. The characteristics of Poisson process lead to the following Proposition:

\textbf{Proposition 1}: With [A1] and [A2], $\text{Pr}(t_{e}^\text{c} = t_{e}^\text{a}~\&~t_{e}^\text{d}-t_{e}^\text{c}\geq T) = \text{Pr}(y^\text{ev}\leq y^\text{cs}),\quad  y^\text{ev}\sim Poisson(T\lambda), \forall e$.

Proposition 1 is intuitive: if the number of PEV arrivals in a duration of $T$ is larger than the number of charging spots, i.e., $y^\text{ev}> y^\text{cs}$, there will be at least one PEV that can not get charged for $T$ units of time with [A2]. A diagram of the Poisson arrivals of homogeneous PEVs is shown in Fig.~\ref{fig:HomoDR}.
The rigorous proof of Proposition 1 is given in the appendix.

\begin{figure}
	\centering
	\tikzset{>=latex}
	\begin{tikzpicture}[scale=0.9]
	\begin{footnotesize}
	\draw[->, thick] (-0.5,0) -- (6.75,0)node[right=1.5pt] {$t$};
	\fill[blue!40!white] (0,0) rectangle (4.75,0.2);
	\fill[red!40!white] (1,0) rectangle (5.75,-0.2);
	
	\node(1)[anchor=south west,inner sep=0] at (-0.75,-1.0-0.1) 
	{\includegraphics[width=0.05\textwidth]{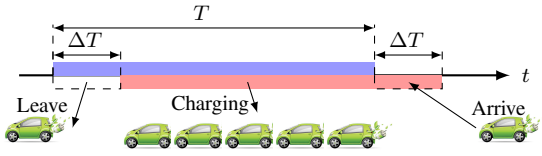}};
	\node(2)[anchor=south west,inner sep=0] at (1.0,-1.2) 
	{\includegraphics[width=0.05\textwidth]{EV.eps}};
	\node(3)[anchor=south west,inner sep=0] at (1.75,-1.2) 
	{\includegraphics[width=0.05\textwidth]{EV.eps}};
	\node(4)[anchor=south west,inner sep=0] at (2.5,-1.2) 
	{\includegraphics[width=0.05\textwidth]{EV.eps}};
	\node(5)[anchor=south west,inner sep=0] at (3.25,-1.2) 
	{\includegraphics[width=0.05\textwidth]{EV.eps}};
	\node(6)[anchor=south west,inner sep=0] at (4.0,-1.2) 
	{\includegraphics[width=0.05\textwidth]{EV.eps}};
	\node(7)[anchor=south west,inner sep=0] at (6.25,-1.0-0.1) 
	{\includegraphics[width=0.05\textwidth]{EV.eps}};
	
	\draw[-,dashed] (0,0) -- (0,0.45);	\draw[-,dashed] (1,0) -- (1,0.45); \draw[<->] (0,0.3) -- (1,0.3)node[above=0pt, xshift=-4.75mm] {$\Delta T$};
	
	\draw[-,dashed] (0+4.75,0) -- (0+4.75,0.45);	\draw[-,dashed] (1+4.75,0) -- (1+4.75,0.45); \draw[<->] (0+4.75,0.3) -- (1+4.75,0.3)node[above=0pt, xshift=-4.75mm] {$\Delta T$};
	
	\draw[-,dashed] (0,0) -- (0,0.75);	\draw[-,dashed] (4.75,0) -- (4.75,0.75); \draw[<->] (0,0.7) -- (4.75,0.7)node[above=0pt, xshift=-23mm] {$T$};
	
	\draw[-, dashed] (0,0) -- (0,-0.2) -- (1, -0.2) -- (1, 0);
	\draw[-, dashed] (4.75,0) -- (4.75,-0.2) -- (5.75, -0.2) -- (5.75, 0);
	
	\draw[->,bend left] (0.5,-0.1) -- (1.east)  node[above=3pt, xshift=-4mm] {Leave};
	\draw[->,bend left] (7.west) -- (5.25, -0.1)  node[below=4pt, xshift=+12mm] {Arrive};
	\draw[->,bend left] (2.875,-0.1) -- (4.north)  node[left=0pt, xshift=-0mm, yshift=1mm] {Charging};
	
	\end{footnotesize}
	\end{tikzpicture}
	\caption{The Poisson arrivals of PEVs (homogeneous driving ranges).}
	\label{fig:HomoDR}
\end{figure}

Based on Proposition 1, Criterion 1 is equivalent to:

\textbf{Criterion 2} \textit{(Homogeneous PEVs)}:  $\text{Pr}(y^\text{ev}\leq y^\text{cs})\geq \alpha, \quad y^\text{ev}\sim Poisson(T\lambda)$.

The Poisson distribution can be approximated by a Normal distribution, i.e., $y^\text{ev}\sim N(T\lambda, T\lambda)$\cite{SupplyChain_Max2011}.
Then, Criterion 2 is:
\begin{align}
&\int_{-\infty}^{y^\text{cs}}f(y^\text{ev})dy^\text{ev}=F(y^\text{cs})=\Phi (\frac{y^\text{cs}-T\lambda}{\sqrt{T\lambda}})\geq \alpha.\label{eqn:station1}
\end{align}

Thus, we have the number of spots in a station limited by:
\begin{align}
&y^\text{cs} \geq F^{-1}(\alpha) = T \lambda + \Phi^{-1}(\alpha) \sqrt{T \lambda}.\label{eqn:station2}
\end{align}

\textcolor{black}{The right-hand side of the above constraint is the closed-form approximation for function $y^\text{cs}(\alpha)$. Constraint (\ref{eqn:station2}) is the service level model subjected to the service level criterion $\alpha$.}

\subsection{The Service Level Model with Heterogeneous PEVs}
In practice, PEVs on the market usually have heterogeneous driving ranges, which results in different charging behaviors, i.e., a PEV with longer driving range may charge fewer times with longer duration each time than a PEV with shorter driving range. Therefore, one single charging station may have to serve PEVs with heterogeneous service time requirements. \textcolor{black}{To effectively model the heterogeneous PEV charging demands, we divide PEVs into different “types” by their driving ranges, e.g., 200 km, 300 km etc. PEVs with similar driving ranges belong to a same ``type" and have similar charging behaviors. We let $\mathcal{K}$ denote the set of PEV types, and $T_k$ denote the required charging time of type $k$ PEVs.}

To derive the closed-form approximation for the function $y^\text{cs}(\alpha)$ of a fast-charging station serving PEVs with heterogeneous driving ranges, we still need to make assumptions [A1] and [A2]. We let $\lambda_k$ denote the Poisson arrival rate of type $k$ PEVs.
Then, similar to Proposition 1, we also have:

\textbf{Proposition 2}: With [A1] and [A2], $\text{Pr}(t_{e_k}^\text{c} = t_{e_k}^\text{a}~\&~t_{e_{k}}^\text{d}-t_{e_{k}}^\text{c}\geq T_k) = \text{Pr}( y^\text{ev}\leq y^\text{cs}),\quad y^\text{ev}=\sum_{k}{y_k^\text{ev}}, y_k^\text{ev}\sim Poisson(T_k\lambda_k),  \forall e_k, \forall k \in \mathcal{K}$.

The proof of Proposition 2 is similar to that of Proposition 1, which is omitted for brevity. A diagram of the Poisson arrivals of two types of PEVs is shown in Fig.~\ref{fig:HeteDR}.

Based on Proposition 2, Criterion 1 is equivalent to:

\textbf{Criterion 3} \textit{(Heterogeneous PEVs)}:  $\text{Pr}( y^\text{ev}\leq y^\text{cs}) \geq \alpha,\quad y^\text{ev}=\sum_{k}{y_k^\text{ev}}, y_k^\text{ev}\sim Poisson(T_k\lambda_k), \forall k \in \mathcal{K}$.

We approximate each independent Poisson distribution $Poisson(T_k\lambda_k)$ by a Normal distribution, i.e.,  $y_k^\text{ev}\sim N(T_k\lambda_{k},  T_k\lambda_{k})$. Because the sum of different independent Normal distributions is still a Normal distribution, we have $y^\text{ev}\sim N(\sum_{k \in \mathcal{K}}T_k\lambda_{k} , \sum_{k \in \mathcal{K}}T_k\lambda_{k})$. Then, Criterion 3 is:
\begin{align}
&\int_{-\infty}^{y^\text{cs}}f(y^\text{ev})dy^\text{ev}=F(y^\text{cs})=\Phi (\frac{y^\text{cs}-\sum_{k \in \mathcal{K}}T_k\lambda_{k} }{\sqrt{\sum_{k \in \mathcal{K}}T_k\lambda_{k} }})\geq \alpha.\label{eqn:station5}
\end{align}

Thus, we have the number of spots in a station limited by:
\begin{align}
&y^\text{cs} \geq F^{-1}(\alpha) = \sum_{k \in \mathcal{K}}T_k\lambda_{k} + \Phi^{-1}(\alpha) \sqrt{\sum_{k \in \mathcal{K}}T_k\lambda_{k}}.\label{eqn:station6}
\end{align}

\textcolor{black}{The right-hand side of the above constraint is the closed-form approximation of $y^\text{cs}(\alpha)$ for a fast-charging station servicing heterogeneous PEVs. Constraint (\ref{eqn:station6}) is the corresponding service level model subjected to the service level criterion $\alpha$.}

\begin{figure}[!t]
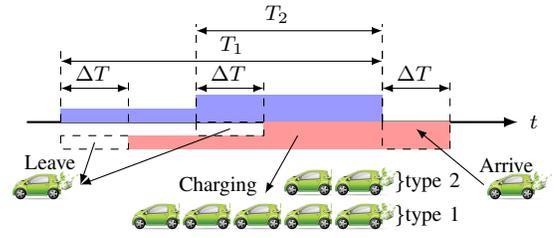

\centering
\tikzset{>=latex}
\begin{tikzpicture}[scale=0.9]
\begin{footnotesize}
\draw[->, thick] (-0.5,0) -- (6.75,0)node[right=1.5pt] {$t$};
\fill[blue!40!white] (2.0,0.2) rectangle (4.75,0.4);
\fill[blue!40!white] (0,0) rectangle (4.75,0.2);
\fill[red!40!white] (3,0) rectangle (5.75,-0.2);
\fill[red!40!white] (1,-0.2) rectangle (5.75,-0.4);

\node(1)[anchor=south west,inner sep=0] at (-0.75,-1.0-0.25) 
{\includegraphics[width=0.05\textwidth]{EV.eps}};
\node(2)[anchor=south west,inner sep=0] at (1.0,-1.7) 
{\includegraphics[width=0.05\textwidth]{EV.eps}};
\node(3)[anchor=south west,inner sep=0] at (1.75,-1.7) 
{\includegraphics[width=0.05\textwidth]{EV.eps}};
\node(3)[anchor=south west,inner sep=0] at (1.75,-1.7) 
{\includegraphics[width=0.05\textwidth]{EV.eps}};
\node(4)[anchor=south west,inner sep=0] at (2.5,-1.7) 
{\includegraphics[width=0.05\textwidth]{EV.eps}};
\node(5)[anchor=south west,inner sep=0] at (3.25,-1.7) 
{\includegraphics[width=0.05\textwidth]{EV.eps}};
\node(6)[anchor=south west,inner sep=0] at (4.0,-1.7) 
{\includegraphics[width=0.05\textwidth]{EV.eps}};

\node(7)[anchor=south west,inner sep=0] at (3.25,-1.2) 
{\includegraphics[width=0.05\textwidth]{EV.eps}};
\node(8)[anchor=south west,inner sep=0] at (4.0,-1.2) 
{\includegraphics[width=0.05\textwidth]{EV.eps}};

\node(9)[anchor=south west,inner sep=0] at (6.25,-1.0-0.25) 
{\includegraphics[width=0.05\textwidth]{EV.eps}};

\draw[-,dashed] (0,0) -- (0,0.65);	\draw[-,dashed] (1,0) -- (1,0.65); \draw[<->] (0,0.5) -- (1,0.5)node[above=0pt, xshift=-4.75mm] {$\Delta T$};

\draw[-,dashed] (2,0) -- (2,0.65);	\draw[-,dashed] (3,0) -- (3,0.65); \draw[<->] (2,0.5) -- (3,0.5)node[above=0pt, xshift=-4.75mm] {$\Delta T$};

\draw[-,dashed] (0+4.75,0) -- (0+4.75,0.65);	\draw[-,dashed] (1+4.75,0) -- (1+4.75,0.65); \draw[<->] (0+4.75,0.5) -- (1+4.75,0.5)node[above=0pt, xshift=-4.75mm] {$\Delta T$};

\draw[-,dashed] (0,0) -- (0,1.05);	\draw[-,dashed] (4.75,0) -- (4.75,1.05); \draw[<->] (0,0.9) -- (4.75,0.9)node[above=0pt, xshift=-20mm] {$T_1$};

\draw[-,dashed] (2,0) -- (2,1.45);	\draw[-,dashed] (4.75,0) -- (4.75,1.45); \draw[<->] (2,1.35) -- (4.75,1.35) node[above=0pt, xshift=-14mm] {$T_2$};

\draw[-, dashed] (0,-0.2) -- (0,-0.4) -- (1, -0.4) -- (1, -0.2) -- (0,-0.2);
\draw[-, dashed] (2,0)--(2,-0.2) -- (3,-0.2) -- (3, -0.0) -- (2, -0.0);
\draw[-, dashed] (4.75,0) -- (4.75,-0.4) -- (5.75, -0.4) -- (5.75, 0);

\draw[->,bend left] (2.5,-0.1) -- (1.east)  node[above=3pt, xshift=-4mm] {};
\draw[->,bend left] (0.5,-0.3) -- (1.east)  node[above=3pt, xshift=-4mm] {Leave};
\draw[->,bend left] (9.west) -- (5.25, -0.1)  node[below=8pt, xshift=+12mm] {Arrive};

\draw[->,bend left] (3.5,-0.2) -- (4.north)  node[left=0pt, xshift=-0mm, yshift=1mm] {Charging};

\draw[-] (6.east) -- (6.east)  node[right=0pt,xshift=-2mm] {\}type 1};
\draw[-] (8.east) -- (8.east)  node[right=0pt,xshift=-2mm] {\}type 2};

\end{footnotesize}
\end{tikzpicture}
\caption{The Poisson arrivals of PEVs (heterogeneous driving ranges).}
\label{fig:HeteDR}
\vspace{-2mm}
\end{figure}

\subsection{Second Order Cone Reformulation of the Model}
\textcolor{black}{Though constraint (\ref{eqn:station6}) is in a simple closed-form, it is non-convex and intractable. In this subsection, we show how to reformulate it into a mixed-integer SOCP. }

For a transportation network with multiple paths, the type $k$ PEV traffic flow (with Poisson arrival rate $\lambda_{q,k}$) to be charged at a node $i$ is given by:
\begin{align}
&\lambda_{i,k}=\sum_{q\in \mathcal{Q}_i}{\lambda_{q,k} \gamma_{q,i,k}},\label{eqn:station7}
\end{align}
where, the charge choice variable $\gamma_{q,i,k}$ is binary; $\gamma_{q,i,k}=1$, if PEVs get charged; $\gamma_{q,i,k}=0$, otherwise.

Thus, the closed-form service level model (\ref{eqn:station6}) for a station at node $i$ serving PEVs with heterogeneous driving ranges is:
\begin{align}
\begin{split}
y_i^\text{cs}\geq & \sum_{q\in \mathcal{Q}_i}\sum_{k \in \mathcal{K}}{T_k\lambda_{q,k} \gamma_{q,i,k}} +\Phi^{-1}(\alpha) \sqrt{\sum_{q\in \mathcal{Q}_i}\sum_{k \in \mathcal{K}}{T_k\lambda_{q,k} \gamma_{q,i,k}} }.
\end{split}\label{eqn:station8}
\end{align}

Because $\gamma_{q,i,k}=\gamma_{q,i,k}^2$ holds, constraint (\ref{eqn:station8}) is equivalent to the following one\cite{Plan_T_Mak2013}:
\begin{align}
y_i^\text{cs}&\geq  \sum_{q\in \mathcal{Q}_i}\sum_{k \in \mathcal{K}}{T_k\lambda_{q,k} \gamma_{q,i,k}}+\Phi^{-1}(\alpha) \sqrt{\sum_{q\in \mathcal{Q}_i}\sum_{k \in \mathcal{K}}{T_k\lambda_{q,k} \gamma_{q,i,k}^2} },\label{singlestation3}
\end{align}
which is a mixed-integer SOCP and can be efficiently solved by the branch-and-cut method in an off-the-shelf commercial solver, such as CPLEX\cite{cplex_SOCP}.


\textbf{Remark 1} Given a service level criterion $\alpha$, constraint (\ref{singlestation3}) provides the minimum required number of spots in a station. The first term in (\ref{singlestation3}) is the required number of spots to satisfy the expected charging demands and is proportional to the Poisson arrival rate. The second term corresponds to the extra spots to satisfy any demand in excess of the mean and can be viewed as the ``safety stock."\footnote{The ``safety stock'' reflects the ``scale effect'' of a station's service ability, i.e., the average required spots per demand decreases as the demands increase.} In practice, high $\alpha$ leads to more investments and ensures better service quality of the charging stations.\footnote{Note that  $\alpha$ is usually required to be above $50\%$ so that the second term (safety stock) in (\ref{eqn:station8}) is above 0\cite{SupplyChain_Max2011}.}

\textbf{Remark 2} The service level model based on assumption [A2] guarantees a lower bound for the future charging station's service quality. 
If the operator of the designed station follows the \textit{first-in first-out rule} in [A2], the service quality is exactly what we have designed.  If the operator smartly manages the station, e.g., allowing waiting or letting the PEV with the highest $SoC$ leave first, the service quality can be higher.

\textcolor{black}{\textbf{Remark 3} The service level criterion, i.e., $\alpha$, provides a meaningful and intuitive service quality criterion for future operations. In practice, waiting is usually allowed; then, $\alpha$ is approximately equal to the probability that a PEV can get instantly charged right after arriving at the station; while $1-\alpha$ is the probability that the PEV has to wait. The probability of waiting in a public charging station is an important service quality criterion in practice and can be easily controlled by $1-\alpha$ in the proposed model. We conduct numerical experiments to validate the performance of this model in Section \ref{sec:model}.}

\section{Modified Capacitated Flow Refueling Location Model Based on Sub-paths}\label{sec:CFRLM}
{This Section describes the transportation network model, i.e., the CFRLM\_SP. It models the PEVs' driving range constraints on the transportation network, which provides the feasible set of the charging station locations in the planning. }
	
Note that the proposed model is based on the origin one first introduced in \cite{Plan_T_Mak2013}. We modify it to make it more accurate and consider time-varying OD traffic flow.

\subsection{PEV Driving Range Constraint Based on Sub-paths}\label{sec:DRlogic}
We assume that the highway networks for PEVs are operated by the following rules:
\begin{enumerate}
	\item Entrance rule: The PEVs should enter the transportation network with battery $SoC$s enough to travel a distance of $D_a$. Because PEVs have limited driving ranges, this rule guarantees the PEVs to be able to reach a charging station without running out of energy.
	\item Exit rule: The PEVs should leave the transportation network with sufficiently high battery $SoC$s enough to travel another distance of $D_d$. This threshold guarantees the service quality of the charging network and ensures that the PEVs leave the network with sufficient energy to arrive at their final destinations.
\end{enumerate}

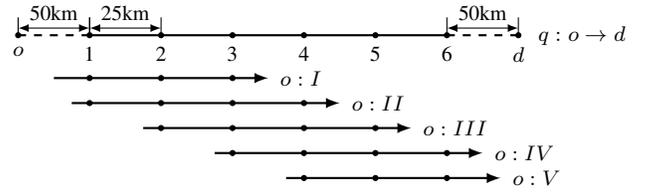
\begin{figure}[!t]
	\centering
		\tikzset{>=latex}
		\begin{tikzpicture}[scale=0.95]
		\begin{footnotesize}
		\draw[-,dashed, thick] (-0.,0) -- (1,0)node[right=1.5pt] {};
		\draw[-,dashed, thick] (6,0) -- (7,0)node[right=4.5pt] {$q:o\rightarrow d$};
		\draw[-, thick] (1,0) -- (6,0)node[right=1.5pt] {};
		
		\filldraw [black] (0,0) circle (1pt) node(1)[below=1.5pt]{$o$};
		\filldraw [black] (1,0) circle (1pt) node(2)[below=1.5pt]{1};
		\filldraw [black] (2,0) circle (1pt) node(3)[below=1.5pt]{2};
		\filldraw [black] (3,0) circle (1pt) node(4)[below=1.5pt]{3};
		\filldraw [black] (4,0) circle (1pt) node(5)[below=1.5pt]{4};
		\filldraw [black] (5,0) circle (1pt) node(6)[below=1.5pt]{5};
		\filldraw [black] (6,0) circle (1pt) node(7)[below=1.5pt]{6};
		\filldraw [black] (7,0) circle (1pt) node(8)[below=1.5pt]{$d$};
		
		\draw[->, thick] (1-0.5,-0.6) -- (3.5,-0.6)node[right=1.5pt] {$o:I$};
		\filldraw [black] (1,-0.6) circle (1pt);
		\filldraw [black] (2,-0.6) circle (1pt);
		\filldraw [black] (3,-0.6) circle (1pt);
		
		\draw[->, thick] (0.75,-0.95) -- (4.5,-0.95)node[right=1.5pt] {$o:II$};
		\filldraw [black] (1,-0.95) circle (1pt);
		\filldraw [black] (2,-0.95) circle (1pt);
		\filldraw [black] (3,-0.95) circle (1pt);
		\filldraw [black] (4,-0.95) circle (1pt);
		
		\draw[->, thick] (1.75,-1.3) -- (5.5,-1.3)node[right=1.5pt] {$o:III$};
		\filldraw [black] (2,-1.3) circle (1pt);
		\filldraw [black] (3,-1.3) circle (1pt);
		\filldraw [black] (4,-1.3) circle (1pt);
		\filldraw [black] (5,-1.3) circle (1pt);
		
		\draw[->, thick] (2.75,-1.65) -- (6.5,-1.65)node[right=1.5pt] {$o:IV$};
		\filldraw [black] (3,-1.65) circle (1pt);
		\filldraw [black] (4,-1.65) circle (1pt);
		\filldraw [black] (5,-1.65) circle (1pt);
		\filldraw [black] (6,-1.65) circle (1pt);
		
		\draw[->, thick] (3.75,-2) -- (6.75,-2)node[right=1.5pt] {$o:V$};
		\filldraw [black] (4,-2) circle (1pt);
		\filldraw [black] (5,-2) circle (1pt);
		\filldraw [black] (6,-2) circle (1pt);
		
		\draw[-] (0,0) -- (0,0.25);
		\draw[-] (1,0) -- (1,0.25);
		\draw[<->] (0,0.1) -- (1,0.1) node[above, xshift=-4.75mm] {50km};
		\draw[-] (6,0) -- (6,0.25);
		\draw[-] (7,0) -- (7,0.25);
		\draw[<->] (6,0.1) -- (7,0.1) node[above, xshift=-4.75mm] {50km};
		\draw[-] (1,0) -- (1,0.25);
		\draw[-] (2,0) -- (2,0.25);
		\draw[<->] (1,0.1) -- (2,0.1) node[above, xshift=-4.75mm] {25km};
		\end{footnotesize}
		\end{tikzpicture}
	\caption{Driving range logic based on sub-path (100 km driving range).}
	\label{fig:sub-path}
\end{figure}

Based on the above operation rules, we explain the driving range logic by Fig.~\ref{fig:sub-path}. A PEV with a driving range of 100 km arrives at node 1 with $D_a=50$ km (which means the PEV has already traveled 50 km before arriving at node 1) and needs to depart at node 6 with $D_d=50$ km. We add pseudo nodes $o$ and $d$ to denote the original node and destination node respectively and let $d_{o,1}=50$ km and $d_{6,d}=50$ km. Then, the problem becomes that a PEV with its battery fully charged leaves at node $o$ and needs to arrive at node $d$ without running out of energy on the road. The travel trajectory of the PEV, i.e., \{$o$, 1, 2, 3, 4, 5, 6, $d$\}, is called a path, i.e., $q$, and a segment of path $q$ is its sub-path. The real nodes on path $q$, i.e., \{1, 2, 3, 4, 5, 6\}, are the candidate locations for charging stations. The \textit{driving range constraint} for a PEV on path $q$ is that any sub-path in $q$ with a distance longer than the PEV's driving range, i.e., 100 km, should cover at least one charging station so that the PEV can travel through path $q$ with adequate charging services. In Fig.~\ref{fig:sub-path}, the set of sub-paths is \{$I, II, III, IV, V$\}. Thus, on every sub-path in \{$I, II, III, IV, V$\}, there should be at least one charging station located on one of its covered nodes. For the case in Fig.~\ref{fig:sub-path}, two stations are required and the candidate locations may be any of \{$1,4$\}, \{$2,4$\}, \{$2,5$\}, \{$3,4$\}, \{$3,5$\}, \{$3,6$\}.

The CFRLM\_SP in \cite{Plan_T_Mak2013} assumed that the PEVs are fully charged before arriving at the highway network, which may not hold in practice. Besides, to consider round trip, \cite{Plan_T_Mak2013} defined the sub-path length to be half of the PEVs' driving ranges, which would make the planning result very conservative. By contrast,we define the sub-path length to be equal to the PEVs' driving ranges to enhance modeling accuracy. Because the PEV traffic flow with round trips should follow the operation rules described above in both the departure and the return trips, we can model both trips separately.


\vspace{-3mm}
\subsection{Capacitated Flow Refueling Location Model}
Based on the driving range constraint described in Section \ref{sec:DRlogic}, the CFRLM\_SP considering time-varying OD traffic flows can be formulated as follows: 
\begin{align}
&\text{Service ability constraint~(\ref{singlestation3})}, \quad \forall i\in \mathcal{I}, \forall \omega \in \Omega , \forall t, \label{eqn:1_1}\\
&\sum_{i\in \mathcal{I}_o}\gamma_{q,i,k} \geq 1, \qquad \forall o\in \mathcal{O}_{q,k}, \forall q\in \mathcal{Q}, \forall k \in \mathcal{K}, \label{eqn:1_2}\\
&x_i \geq \gamma_{q,i,k}, \qquad \forall q\in \mathcal{Q}, i\in \mathcal{I}, \forall k \in \mathcal{K},\label{eqn:1_4}\\
&0 \leq y_i^\text{cs}\leq \overline{y_{i}}x_{i}, \qquad \forall i\in \mathcal{I}.\label{eqn:1_5}
\end{align}
Each station's service ability is constrained by (\ref{eqn:1_1}). To consider time-varying OD traffic flows, the service ability constraint (\ref{eqn:1_1}) should be satisfied for every hour in every scenario.\footnote{We utilize the hourly average traffic flow, i.e., $\lambda_{q,i,k,\omega,t}$, in this constraint. }
\textcolor{black}{Equation (\ref{eqn:1_2}) ensures that the PEVs get charged for at least once in each sub-path. This is the formulation of the driving rang constraint introduced in Section \ref{sec:DRlogic}.}
 Equation (\ref{eqn:1_4}) constrains that the PEVs can only get charged at the nodes with charging stations. Equation (\ref{eqn:1_5}) upper-bounds the number of charging spots, if there is a charging station at that node.

\subsection{Extra Constraints for CFRLM\_SP}\label{sec:cut}

\begin{figure}
	\centering
	\tikzset{>=latex}
	\begin{tikzpicture}[scale=1]
	\begin{footnotesize}
	\draw[->, thick] (-0.25,0) -- (6.75,0)node[right=1.5pt] {$q_1$};
	
	\filldraw [black] (0,0) circle (1pt) node(1)[below=1.5pt]{1};
	\filldraw [black] (1,0) circle (1pt) node(2)[below=1.5pt]{2};
	\filldraw [black] (2,0) circle (1pt) node(3)[below=1.5pt]{3};
	\filldraw [black] (3,0) circle (1pt) node(4)[below=1.5pt]{4};
	\filldraw [black] (4,0) circle (1pt) node(5)[below=1.5pt]{5};
	\filldraw [black] (5,0) circle (1pt) node(6)[below=1.5pt]{6};
	\filldraw [black] (6,0) circle (1pt) node(7)[below=1.5pt]{7};
	
	\filldraw [black] (5,-0.6) circle (1pt)
	node(8)[below=1.5pt]{8};
	\filldraw [black] (6,-0.6) circle (1pt) node(9)[below=1.5pt]{9};
	\draw[->, thick] (4,0) -- (5,-0.6) -- (6.75,-0.6) node[right=1.5pt] {$q_2$};
	
	\draw[-] (0,0) -- (1,0) node[above=1.5pt, xshift=-4.75mm] {25km};
	\draw[-] (0,0) -- (0,0.25);
	\draw[-] (1,0) -- (1,0.25);
	\draw[<->] (0,0.1) -- (1,0.1);
	
	\end{footnotesize}
	\end{tikzpicture}
	\vspace{-2mm}
	\caption{Two paths, i.e., $q_1$ and $q_2$, with a same arrival node 1 and an identical sub-path, i.e., \{1, 2, 3, 4, 5\}. }
	\label{fig:twopath}
\end{figure}
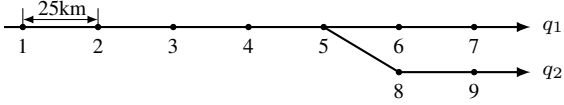

\textcolor{black}{In the CFRLM\_SP, the scale of the binary variables increase linearly with the scale of the transportation network. In practice, the transportation network may be complicated so that the optimization model may be intractable. To decrease the scale of the problem, we make the following mild assumption: }

\textbf{[A3]} The PEVs with the same arrival node traveling on identical sup-paths have the same charge choices before they separate with each other at the end of the identical sub-paths. 

Take Fig.~\ref{fig:twopath} as an example. PEVs traveling from node 1 to node 7 (path $q_1$) have the same charge choices between node 1 and 5 with those PEVs traveling from node 1 to node 9 (path $q_2$). This assumption actually has practical meaning that the PEVs on path $q_1$ and $q_2$ will have the same traveling experiences before they depart with each other at node 5. Therefore, they tend to have the same charge choices when they visit a charging station. 

With [A3], the scale of binary variables, i.e., $\gamma_{q,i,k}$, decreases significantly so that the computational efficiency is enhanced. Note that the reduction of the scale of $\gamma_{q,i,k}$ depends on the structure of the transportation network. Take a scenario with two paths, e.g., Fig.~\ref{fig:twopath}, for example. If the identical sub-path has $o_0$ candidate locations, and the separate sub-paths respectively have $o_1$ and $o_2$ variables, then, the scale of $\gamma_{q,i,k}$ is reduced from $o_0\times 2+o_1+o_2$ to $o_0+o_1+o_2$ with [A3]. Apparently, if different paths share longer identical sub-paths, the scale of $\gamma_{q,i,k}$ will reduce more significantly. When the lengths of the separate sub-paths, i.e., $o_1$ and $o_2$, are very short compared to that of the identical path, i.e., $o_0$, the two paths can be regarded as a single path.

\section{Charging Station Planning Model}
{As introduced earlier, the service ability model in Section \ref{sec:SRM} determines  the sizes of the charging stations given the demands; and the transportation network model (the CFRLM\_SP) in Section \ref{sec:CFRLM} defines the feasible set of the PEV charging sites. This section formulates the complete model for siting and sizing PEV fast-charging stations subjected to the aforementioned two models. Furthermore, we introduce the AC power flow model to describe the power network constraints' impact on the planning.  
}
\subsection{Objective}
Considering that the traffic flow and traditional base loads are uncertain over the target planing horizon, a set of finite potential future scenarios ($\Omega$) are forecasted. Then a  two-stage stochastic programming model is adopted to plan fast-charging stations. 
The objective is formulated as follows:
\begin{align}
&Obj=\zeta\sum_{i\in \mathcal{I}}\left(c_{1,i}x_i+c_{2,i}y_i^\text{cs} +c_{3,i}l_{i}\overline{P_{i}^\text{ev}}+c_{4,i}P_i^\text{sub}\right)\notag\\
&+365\sum_{\omega\in \Omega}\sum_{t} \pi_\omega\left( c_\text{e} p_{0,\omega,t} \Delta t + c_\text{p} \sum_{i\in \mathcal{I}}{p_{\text{un},i,\omega,t}}\Delta t \right),\label{obj2}
\end{align}
where:
\begin{align}
&\overline{{P}_{i}^\text{ev}}=p^\text{spot} y_i^\text{cs}, \qquad \forall i\in \mathcal{I},\label{eqn:2-1}\\
&P_i^\text{sub}=\max(0,\overline{{P}_{i}^\text{ev}}-P_{i,0}^\text{sub}), \qquad \forall i\in \mathcal{I}.\label{eqn:2-2}
\end{align}

The first two terms in the first line of (\ref{obj2}) represent the fixed cost of building charging stations and the variable building cost in proportion with the number of charging spots. The last two terms in the first line of  (\ref{obj2}) together account for power distribution network upgrade costs, which include the costs for distribution lines and the costs for substation capacity expansion. 
The first term in the second line is the annual expected energy purchase costs of the whole system and the second term is the penalty for unsatisfied charging demands.
The maximum charging power in each station, i.e., $\overline{P_{i}^\text{ev}}$, is calculated by (\ref{eqn:2-1}). The corresponding substation capacity expansion is calculated by (\ref{eqn:2-2}).


\subsection{Constraints}
\subsubsection{Transportation Network Constraints} The model should satisfy the constraints of CFRLM\_SP, i.e., (\ref{eqn:1_1})--(\ref{eqn:1_5}).

\subsubsection{Power Network Constraints}\label{sec:EC}
\textcolor{black}{The branch currents and nodal voltages of the distribution network must satisfy the AC power flow constraints. In this paper, the SOCP relaxation of AC power flow model \cite{OPF_tree_Exactness_Gan2015} is adopted, as follows:}
\begin{align}
&\forall m\in\mathcal{M}, \forall \left(m,n\right)\in \mathcal{L}, \forall \omega \in \Omega, \forall t: \notag\\
&S_{mn,\omega,t}=s_{m,\omega,t}+\sum_{h\in \mathcal{M}_{\rightarrow m}}{\left(S_{hm,\omega,t}-z_{hm}l_{hm,\omega,t}\right)}, \label{eqn:opf1}\\
&0=s_{0,\omega,t}+\sum_{h\in \mathcal{M}_{\rightarrow 0}}{\left(S_{h0,\omega,t}-z_{h0}l_{h0,\omega,t}\right)}, \label{eqn:opf2}\\
&v_{m,\omega,t}-v_{n,\omega,t}=2\text{Re}(z_{mn}^{*}S_{mn,\omega,t})-|z_{mn}|^2l_{mn,\omega,t}, \label{eqn:opf3}\\
&|S_{mn,\omega,t}|^2 \leq l_{mn,\omega,t}v_{m,\omega,t}, \label{eqn:opf4}\\
&s_{m,\omega,t}=-s_{m,\omega,t}^\text{ev}-s_{m,\omega,t}^\text{b}.
\end{align}
Note that, in this paper, we assume that the coupled power system is a radial high-voltage distribution network. Considering that the system only have unidirectional power flow, the SOCP relaxation is exact \cite{OPF_tree_Exactness_Gan2015}. 
For scenarios when higher voltage level transmission systems are also covered in the targeted planning area, using the SOCP relaxation for the whole system may not provide a feasible solution for the problem. In that case, we can adopt the linear direct current (DC) power flow to model the transmission networks\footnote{\textcolor{black}{In transmission systems, the line resistances are negligible compared to the reactances, the per-unit voltage amplitudes of different nodes are approximately equal to 1.0, and the voltage angle differences between neighboring nodes are small. As a result, the DC power flow model is an accurate approximation for the AC power flow model\cite{DCPowerFlow_stott2009dc}.}}; meanwhile, we can still apply the SOCP relaxation to model the AC power flow of the radial distribution networks. This modification is easy to implement and will not impact the characteristics of the planning model. For brevity, we omitted the corresponding formulations in this paper. 

The distribution line currents and nodal voltages must not violate their permitted ranges:
\begin{align}
&l_{mn,\omega,t}\leq |\overline{I_{mn}}|^2, \qquad \forall \left(m,n\right)\in \mathcal{L}, \forall \omega \in \Omega, \forall t,\label{eqn:current}\\
&|\underline{V_{m}}|^2\leq v_{m,\omega,t}\leq |\overline{V_{m}}|^2, \qquad \forall m\in\mathcal{M}, \forall \omega \in \Omega, \forall t.
\end{align}

\subsubsection{Coupled Constraints of Transportation \& Power Networks}
The PEV charging load at each distribution bus is calculated as follows:
\begin{align}
&s_{m,\omega,t}^\text{ev}=p_{m,\omega,t}^\text{ev} + j \tan\theta  p_{m,\omega,t}^\text{ev}, ~\forall m \in\mathcal{M}, \forall \omega \in \Omega, \forall t,\\
& p_{m,\omega,t}^\text{ev}= \sum_{i\in \mathcal{I}_{m}}p_{i,\omega,t}^\text{ev}, \qquad \forall m \in \mathcal{M}, \forall \omega \in \Omega, \forall t,
\end{align}
where, $\theta$ is the phase angle between the PEV charging voltage and current; $\cos \theta$ is the power factor of the PEV charging load. 
Considering that battery chargers usually have high power factors (close to 1.0, see \cite{Powerfactor_zhang2012unit,Powerfactor_hou2013comparison}), we can approximately assume $\cos \theta =1$ so that $\tan \theta =0$.

The average PEV charging power at each transportation node is calculated as follows:
\begin{align}
&p_{i,\omega,t}^\text{ev} + p_{\text{un},i,\omega,t}^\text{ev}= p^\text{spot}\sum_{q\in \mathcal{Q}_i}\sum_{k \in \mathcal{K}}{T_k\lambda_{q,i,k,\omega,t} \gamma_{q,i,k}},\notag \\
& \qquad \qquad \qquad \qquad ~ \forall i \in \mathcal{I}, \forall \omega \in \Omega, \forall t.\label{eqn:coupled2}
\end{align}

Note that when the PEV traffic is low, $P_{\text{un},i,\omega,t}^\text{ev}=0$ and the PEV charging power $p_{i,\omega,t}^\text{ev}$ is proportional to the average traffic flow that required charging services. On the other hand, when the traffic flow grows beyond the system's service ability, some charging demands are not fulfilled and $P_{\text{un},i,\omega,t}^\text{ev}> 0$, which reveals the power network's influence on the planning.

The base loads are required to be satisfied in the model, i.e., $s_{m,\omega,t}^\text{b}=q_{m,\omega,t}^\text{b}+jq_{m,\omega,t}^\text{b}, \forall m \in \mathcal{M}, \forall \omega \in \Omega, \forall t$.

The above formulated model, i.e., objective (\ref{obj2}) subjected to constraints (\ref{eqn:1_1})--(\ref{eqn:1_5}) and (\ref{eqn:2-1})--(\ref{eqn:coupled2}), is an MISOCP. It can be directly solved by the branch-and-cut method in a commercial solver such as CPLEX\cite{cplex_SOCP}. Because the charging demands can be unsatisfied when the system operation constraints, e.g., AC power flow, are binding, the model is always feasible.

\begin{figure}[!t]
	\centering
	\vspace{-5mm}
	\includegraphics[width=0.7\columnwidth]{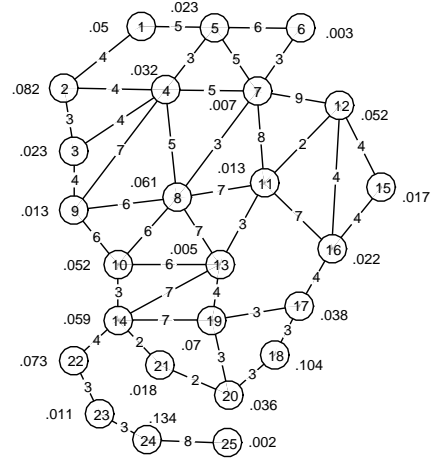}
	\vspace{-9mm}
	\caption{A 25-node transportation network used for the case study\cite{Zhang_PlanFRLMTSG2016}. The number in each circle is the node ID. The number on each line represents the distance between the corresponding two nodes and the per-unit distance is 10 km. The decimal next to each node is its weight, i.e., $W$, which represents its traffic flow gravitation\cite{Zhang_PlanFRLMTSG2016}. To enhance network granularity, we add extra auxiliary nodes on the long line segments so that the longest distance between any two nodes is 20 km. As a result, the modified network has 93 nodes.}
	\label{fig:trans25}
\end{figure}
\begin{figure}[!t]
	\centering
	\includegraphics[width=0.85\columnwidth]{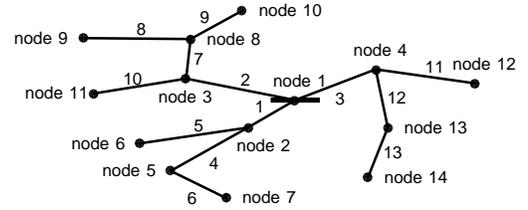}
	\vspace{-6mm}
	\caption{A 110 kV distribution network used for the case study\cite{Zhang_PlanFRLMTSG2016}. Node 1 is connected to a 220 kV/110 kV transformer with 150 MVA capacity. The voltage constraints are $\underline{V_{m}}=0.95$ and $\overline{V_{m}}=1.05$, $\forall m$, in per unit values. The line current limits are conservatively set at 85\% of their rated capacities. The detailed parameters of the distribution network are given in \cite{Supplementary_Material_SOCP_Planning_Model}.}
	\label{fig:distribution14}
\end{figure}
\begin{table}[!t]
	\renewcommand{\arraystretch}{1.1}
	\centering
	\begin{small}
		\caption{Node Coupling Relationship of the Two Networks}
		\begin{tabular}{cccccccc}
			\hline\hline
			Distribution Node ID&01&02&03&04&05&06&07\\
			Transportation Node ID&--&13&08&12&22&14&24\\
			\hline
			Distribution Node ID&08&09&10&11&12&13&14\\
			Transportation Node ID&04&02&05&09&15&17&20\\
			\hline\hline
		\end{tabular}\label{tab:coupling}
	\end{small}
\end{table}

\section{Case Studies}
\subsection{Case Overview and Parameter Settings}
We consider a 25-node highway transportation network (see Fig.~\ref{fig:trans25}) coupled with a 14-node 110 kV high voltage distribution network (see Fig.~\ref{fig:distribution14}) to illustrate the proposed planning method. 
Note that we adopt the power system structure in China as the basis of this case study, where the 110 kV power networks are usually operated radially and categorized as high-voltage distribution systems\cite{ChinaDistribution_Zhong2007}.
The node coupling relationship between the distribution and transportation network is recorded in Table~\ref{tab:coupling}. We assume the transportation nodes not included in Table~\ref{tab:coupling} are connected to the nearest distribution nodes geographically. The gravity spatial interaction model utilized in \cite{Zhang_PlanFRLMTSG2016} was used to generate a daily OD flow structure based on node weights and arc distances. 
Twenty-four representative scenarios, i.e., weekday and weekend of 12 months, of base load profiles and traffic flow profiles are generated based on PG\&E load profiles\cite{Loadprofile_PGE} and NHTS data\cite{NHTS}. 
Due to limited space, the parameters of the distribution network and the details of the generated scenarios are omitted, but can be downloaded in \cite{Supplementary_Material_SOCP_Planning_Model}.

We assume four types of PEVs on road with equal market share, and their driving ranges per charge are respectively 200, 300, 400 and 500 km. We assume the energy consumption of all types of PEV are all 0.14 kWh/km\cite{Plan_Zhang2015}. The rated charging power for each charging spot ($p$) is 44 kW, and the charging efficiency ($\eta$) is 92\%\cite{Plan_Zhang2015}. Consequently, the average service time to recharge the four different types of PEVs with empty batteries is about 42, 63, 84, 105 minutes.  We also assume $D_a=D_d=100$ km for all PEVs and $\overline{y_{i}}=200$. 

The costs of charging station investment $c_{1,i}=\$ 163,000$ and $c_{2,i}=\$ 31,640$.
The distribution line cost $c_{3,i}=120$ \$/(kVA$\cdot$km)\cite{Paramater_Yao2015}. The distance from the PEV charging station to its nearest distribution substation, i.e., $l_{i}$, is assumed to be 10\% of the distance between the PEV charging station and its nearest 110 kV distribution node. The substation expansion cost $c_{4,i}=788$ \$/kVA\cite{Plan_TE_Yao2014}. In practice, the land use and labor costs vary by location. To model this feature across nodes, the per-unit costs, i.e., $c_{1,i}$, $c_{2,i}$ and $c_{4,i}$, at each location $i$ are assumed to be greater than the base values introduced above by $5W_{i}\times100\%$. We assume each original transportation node has 1 MVA surplus substation capacity which can be utilized by charging stations, while the auxiliary nodes have no spare capacity. 
The electricity purchase cost $c_\text{e}=0.094$ \$/kWh \cite{Plan_Zhang2015} and the penalty cost for unsatisfied charging demand $c_\text{p}=10^{3}$ \$/kWh. The service quality $\alpha=80\%$. 

\textcolor{black}{Note that the above parameters are for illustration purposes. In practice, the planner should adopt the actual PEV parameters in current and future markets based on a practical business survey and substitute their own parameter values for the transportation and power networks.}

\begin{table*}[!t]
	\renewcommand{\arraystretch}{1.1}
	\centering
	\begin{small}
		\caption{Benchmark Cases of the Planning}
		\begin{tabular}{cccccc}
			\hline\hline
			\multirow{2}{*}{Case}&\multirow{2}{*}{Driving range} & Extra constraints & \multirow{2}{*}{Electrical constraints} &Distribution system&Traffic flow per day in \\
			& &in Section \ref{sec:cut}& & upgrade cost& the highest traffic scenario\\
			\hline
			1&heterogeneous&Consider&AC power flow&Consider&20000\\
			2&heterogeneous&Consider&AC power flow&Consider&40000\\
			\hline
			3&homogeneous&Consider&AC power flow&Consider&20000\\
			4&heterogeneous&Ignore&AC power flow&Consider&20000\\
			\hline
			5&heterogeneous&Consider&DC power flow&Consider&20000\\
			6&heterogeneous&Consider&DC power flow&Consider&40000\\
			\hline
			7&heterogeneous&Consider&Ignore&Ignore&20000\\
			8&heterogeneous&Consider&Ignore&Ignore&40000\\
			\hline\hline
		\end{tabular}\label{tab:case}
	\end{small}
\end{table*}

{We design eight cases to illustrate the proposed method (see Table~\ref{tab:case}). Case 1 is the basic case utilizing the proposed method. 
In Case 2, the daily PEV traffic flow is twice of that in Case 1. 
In Case 3, we ignore the heterogeneity of the PEV driving ranges. For planning purposes, we conservatively assume all the PEVs are homogeneous and have the shortest driving range, i.e., 200 km, as \cite{Zhang_PlanFRLMTSG2016} suggested.
Case 4 does not consider the extra constraints introduced in Section \ref{sec:cut}.
In Cases 5 and 6, we adopt the DC approximation for the power flow constraints which was also used in reference \cite{Zhang_PlanFRLMTSG2016}.
In Cases 7 and 8, the electrical constraints and the distribution system upgrade costs are ignored at the planning stage.}

We use CPLEX\cite{cplex_SOCP} to solve the optimal PEV charging station planning problem on a laptop with a 12 core Intel Xeon E5-1650 processor and 64 GB RAM. To accelerate the optimization speed, we relaxed the integer variable constraints for the number of spots $y_i^\text{cs}$. The optimization problem stops when the relevant gap decreases below 0.5\%.

\subsection{Planning Results and Analysis}
{The summary of the planning results for the eight cases are given in Table~\ref{tab:result}.
In Cases 3, 5 and 6, the electricity cost, the total cost, and the ratios of unsatisfied PEV load in the parentheses are the direct outputs (solutions) of the optimization models. However, these solutions do not reflect the true operation scenarios. After the investment decisions are obtained, we conduct extra optimizations, which utilize the proposed model in Case 1 but fix the investment decisions, to calculate their actual values. These values are listed outside of the parentheses in Table~\ref{tab:result}. For Case 3, the actual values are calculated considering the heterogeneity of the PEV driving ranges. For Cases 5 and 6, they are calculated by modeling AC power flow. 
Because the electrical constraints in Cases 7 and 8 are ignored at the planning stage, the optimization models will not calculate the grid upgrade cost, the electricity cost and the ratios of unsatisfied PEV load for the planning results. The corresponding values in Table~\ref{tab:result} are also calculated by the model in Case 1 after the planning results, i.e., the number of charging stations and the number of charging spots in each station, are given.
The site and size of each station in case 1 are given in Fig.~\ref{fig:station1} for demonstration.}

{Hereinafter, we will discuss the impact of various factors on the planning results comparing the above eight cases. These factors include the PEV population, the heterogeneity of PEV driving ranges, and the extra constraints introduced in Section~\ref{sec:cut} etc. We will also discuss the modeling accuracy of AC power flow in Section \ref{sec:DC} and the necessity of considering the coupled constraints of transportation and power networks in Section \ref{sec:coupled}.}

\begin{figure}[!t]
	\centering
	\vspace{-6mm}
	\includegraphics[width=0.7\columnwidth]{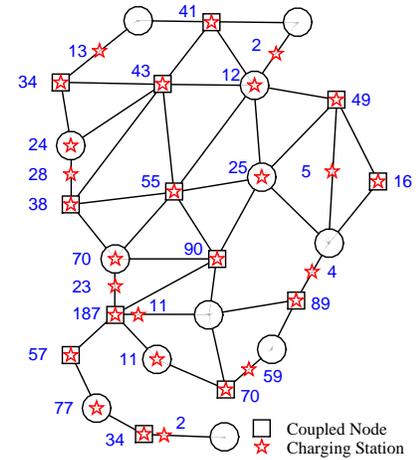}
	\vspace{-8mm}
	\caption{Planning result for Case 1. The integer represents the number of charging spots in the charging station.}
	\label{fig:station1}
	\vspace{-4mm}
\end{figure}

\begin{table*}[!t]	\renewcommand{\arraystretch}{1.1}
	\centering
	\begin{footnotesize}
		\begin{threeparttable}
		\caption{The planning results of different cases}
		\begin{tabular}{cccccccc|cc}
			\hline\hline
			\multirow{2}*{Case}&No. of &No. of  &\multicolumn{4}{c}{Expected annual costs (M\$)} & Unsatisfied & No. of binary  & Solution\\
			&stations&spots&Station Investment & Grid upgrade & Electricity & Total & PEV load (\%) & variables& time (min)\\
			\hline
			1&28&1169&5.13&4.37&38.33&47.83&0&5,761&9.37\\
			2&45&2340&9.91&12.41&49.27&71.59&1.89&5,761&11.10\\
			\hline
			3&46&2722&11.30&14.86&38.07 (51.35)&64.32 (77.50)&0 (5.76)&1,509&8.10\\
			4&19&1017&4.40&2.93&37.43&44.75&0&21,757&92.02\\
			\hline
			5&29&1183&5.20&4.54&38.36 (42.53)&48.10 (52.27)&0 (0.006)&5,761&8.01\\
			6&47&2308&9.85&12.19&48.94 (52.50)&70.97 (74.54)&2.04 (3.55)&5,761&10.35\\
			\hline
			7&21&1160&4.88&6.47&38.45&49.80&0.034&5,761&16.81\\
			8&26&2211&8.99&12.64&48.27&69.89&4.53&5,761&30.67\\
			\hline\hline
		\end{tabular}\label{tab:result}
		\begin{tablenotes}
			\item \textbf{Note:} In Cases 3, 5 and 6, the values in the parentheses are the direct outputs of the optimization models while those outside are the actual values.
		\end{tablenotes}
	\end{threeparttable}
	\end{footnotesize}
\end{table*}
\begin{figure*}[!]
	\centering
	\vspace{-0mm}
	\begin{minipage}[t]{1\columnwidth}
		\centering
		\includegraphics[width=0.9\columnwidth]{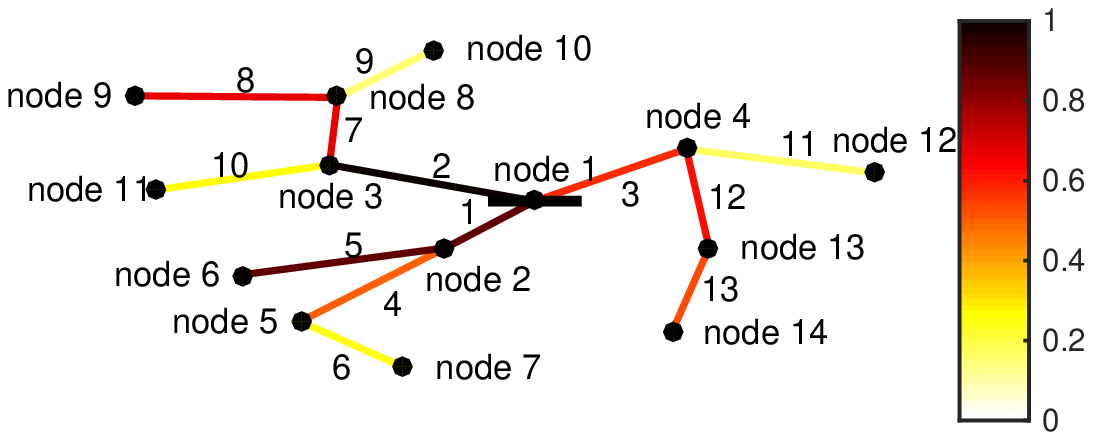}
		\vspace{-3mm}
		\caption{Distribution line congestion level (20000 PEVs/day).}
		\label{fig:grid1}
	\end{minipage}
	\hspace{1mm}
	\begin{minipage}[t]{1\columnwidth}
		\centering
		\includegraphics[width=0.9\columnwidth]{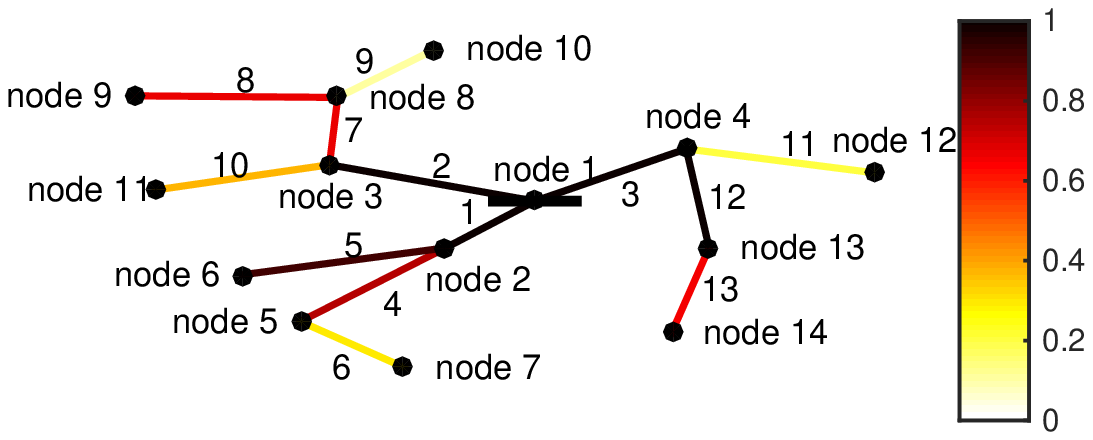}
		\vspace{-3mm}
		\caption{Distribution line congestion level (40000 PEVs/day).}
		\label{fig:grid2}
	\end{minipage}
\end{figure*}

\subsubsection{PEV population}
In the long term, the PEV population is uncertain, some sensitivity analysis is necessary. Compared with case 1, the PEV population in Case 2 is increased by 100\% so that its investment and operation costs both increase significantly. Furthermore, a noteworthy portion of PEV charging load, i.e., 1.89\%, is unsatisfied, which reveals the impact of the power network's constraints. The power flow congestion level, i.e., the ratios of the distribution lines' currents to their capacities, at 12:00 p.m., weekday, July of Case 1 and Case 2 are respectively depicted in Fig.~\ref{fig:grid1} and Fig.~\ref{fig:grid2}. Obviously, the congestion level of Case 2 is much more serious than that of Case 1. Besides, we can observe that the capacity of distribution line 2 is the bottleneck of the system, which provides guidance for future distribution system expansion.

\subsubsection{Heterogeneity of PEV driving range}
When all the PEVs are assumed to have homogeneous driving ranges, the model tends to construct more charging stations and spots, which leads to a very conservative planning result. Compared with Case 1, the total investments in Case 3 is about 100\% higher. The estimated electricity consumption costs also increase significantly. This is because when assuming all PEVs have the same shortest driving range, those with longer driving ranges would charge more times than actually needed. Similarly, the estimated unsatisfied PEV load ratio in Case 3 is very high, i.e., 5.76\%, because of over-estimated demands. However, in practice, since the planning result is very conservative, the PEVs will get sufficiently charged and the spots will be under-utilized. That would be a waste of investments.

\subsubsection{Extra constraints in Section~\ref{sec:cut}} 
Case 4 does not consider the constraints introduced in Section~\ref{sec:cut}, so that the planning results are very aggressive: the investment costs are reduced by 22.8\% and the total costs are reduced by 6.4\%. The physical meaning of these gaps are the revenue that the system can reap by smartly navigating all the PEVs' charge choices so that the utilization of charging stations can be maximized. Note that, in practice, the planner should also consider the smart navigation system's costs in the planning model in Case 4 to make an economic decision, i.e., whether invest the navigation system or not. However, if the smart navigation system is not guaranteed to be implemented and the PEVs may not adopt their optimal (for the whole system) charge choices, this planning strategy may lead to congestion at some stations. To promote the future PEV adoption and guarantee adequate service quality, it is wise to make a conservative infrastructure investment plan.
Regarding the computational efficiency, the scale of binary variables in Case 4 is about four times of that in Case 1. Hence, Case 4's solution time is much longer. In reality, the highway networks may be very complicated and will limit the applicability of the strategy utilized in Case 4.
Therefore, we recommend the proposed strategy with the extra constraints introduced in Section~\ref{sec:cut}.

\subsubsection{Modeling accuracy of AC power flow}\label{sec:DC}
The DC approximation for the power flow model in distribution systems is inaccurate. We can observe that the actual values (utilizing AC power flow) and the estimated values (utilizing DC power flow) in Cases 5 and 6 have apparent differences. The electricity costs in both Cases 5 and 6 are overestimated by about 10\%.  While, in Case 5, though the PEV charging demands can be fully satisfied, we observe that 0.006\% of the demands are unsatisfied in the optimization model. Similarly, the ratio of unsatisfied charging demands in Case 6 is overestimated by about 74\%.
This inaccuracy will lead to sub-optimal planning decisions. 
Compared with Case 1, the actual investment cost and the actual total annual cost in Case 5 increased by 2.5\% and 0.6\%, respectively.
Compared with Case 2, the actual ratio of unsatisfied PEV charging load in Case 6 increased by 7.9\%. Therefore, if possible, it is better to utilize AC power flow in the planning model.

\subsubsection{Necessity of considering coupled transportation and power network constraints}\label{sec:coupled}
When ignoring the power network constraints, i.e., the AC power flow constraints, the unsatisfied PEV charging demands will increase. In Case 1, all the PEV charging demands are satisfied; while, in Case 7, 0.034\% of the PEV charging demands can not be fulfilled. Compared with Case 2, the ratio of unsatisfied PEV charging demands of Case 8 increased significantly by 139.7\%. When the PEV population is large and the power supply capacities are binding in some distribution nodes, the planner prefers to invest more charging stations elsewhere to avoid distribution system congestion. As a result, the numbers of charging stations and spots in Case 2 are much higher than those in Case 8. Besides, the total investment costs in Case 2 are also higher than that in Case 8 (in order to satisfy more charging demands).

When ignoring the distribution system upgrade costs in the planning model, the planner may conduct myopic investment decisions based on limited information. As a result, though the PEV charging station investment costs may be minimized, the system planner has to invest more in distribution system upgrades, which will surpass the savings in PEV charging station investments. This is predominant under low PEV population scenarios when almost all the PEV charging demands can be satisfied (see Case 1 and Case 7). Considering distribution system upgrade costs helps reduce total investment costs.

\subsection{Performance of the Service Level Model}\label{sec:model}
We utilize the proposed model to design the number of charging spots in a fast-charging station under different PEV Poisson arrival parameters, i.e., from 20 to 300 PEVs/h, and different service level criterion, i.e., from 70\% to 90\%. 
Then, we utilize the Monte-Carlo method\cite{MonteCarlo_mooney1997} to simulate the real-time operation including PEV arriving, charging and leaving behaviors of the designed station for 1000 hours.

\begin{figure}[!t]
	\centering
	\includegraphics[width=0.9\columnwidth]{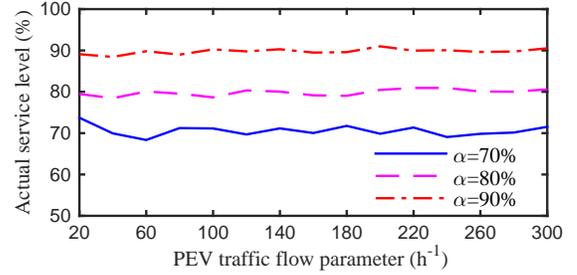}
	\caption{Accuracy of the service level model.}
	\label{fig:servicerate}
\end{figure}
\subsubsection{Accuracy of the service level model} We first assume that the charging stations are operated based on the \textit{first-in-first-out} rule in assumption [A2]. We counted the number of unsatisfied charging demands (those leave before getting charged for their required time units), the corresponding actual service levels are plotted in Fig.~\ref{fig:servicerate}.
The actual service levels in the experiments are very close to the designed values under different PEV Poisson arrival parameters and designed service level criterion. This demonstrates the accuracy of the proposed closed-form approximation for the service level model.

\begin{figure}[!t]
	\centering
	\includegraphics[width=0.9\columnwidth]{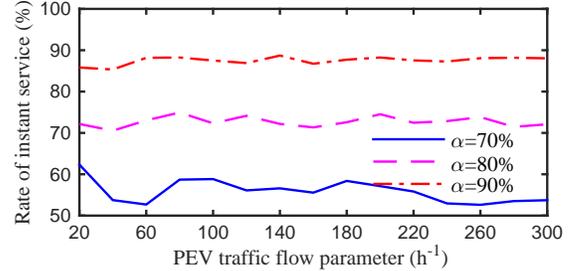}
	\caption{Probabilities that the PEVs get instant charging services.}
	\label{fig:serviceratewaiting}
\end{figure}

\begin{figure}[!t]
	\centering
	\includegraphics[width=0.9\columnwidth]{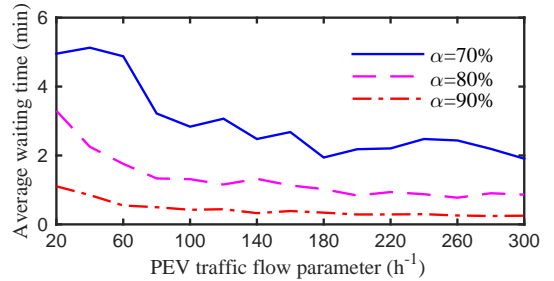}
	\caption{The average waiting time of all the PEVs.}
	\label{fig:waiting}
\end{figure}
\subsubsection{Waiting time analysis in real-time operations} We then assume that the PEVs will wait in the station if all the charging spots are occupied until one PEV has got fully charged and spare a spot for it. We calculated two service quality criterions, i.e., the ratios that the PEVs get instantly charged without waiting (see Fig.~\ref{fig:serviceratewaiting}), and the average waiting time of all the PEVs (see Fig.~\ref{fig:waiting}).

On one hand, the experiment results show that the probabilities that the PEVs get instantly charged without waiting are close to the designed service level criterion. This is especially true when the designed service level criterion is high and few PEVs have to wait. Even though the two values, i.e., the probability of instant services and the designed service level criterion, are not equal, the former is still well controlled by the latter. The results show that the differences between the two values under different designed service level criterion are stable and insensitive to the volumes of the PEV arrivals. That means, we can realize an one-to-one mapping between the two values so that we can use the service level to accurately control the actual probability of waiting.

On the other hand, we can also conclude from Fig.~\ref{fig:waiting} that the average waiting time of the PEVs are also well controlled by the designed service level criterion. When the service level criterion increases, the average waiting time decreases.

Therefore, though the assumption [A2] may not be true in practice, the service level criterion in the proposed model is still meaningful because it provides the planner an intuitive and specific service quality criterion for the designed system.

\begin{table}[!t]
	\renewcommand{\arraystretch}{1.05}
	\centering
	\begin{small}
	\begin{threeparttable}
		\caption{Computationally efficiency of different models}
		\begin{tabular}{ccccc}
			\hline\hline
			\multirow{2}*{Case}&\multirow{2}*{Model}&No. of   & No. of binary& Solution  \\
			&& PEV types&  variables & time (min) \\
			\hline
			1&CFRLM\_SP&4& 5761 & 9.37 \\
			3&CFRLM\_SP&1& 1509 & 8.10\\
			1&CFRLM\_EN&4& 74441 & $\infty$\\
			3&CFRLM\_EN&1& 15198&  120\\
			\hline\hline
		\end{tabular}\label{tab:scale}
		\begin{tablenotes}
		\begin{footnotesize}
			\item \textbf{Note:}  The program in Case 1 utilizing the CFRLM\_EN was out of memory.
			\end{footnotesize}
		\end{tablenotes}
	\end{threeparttable}
	\end{small}
\end{table}

\subsection{Computationally Efficiency Compared with Transportation Model CFRLM\_EN}
In our previous work \cite{Zhang_PlanFRLMTSG2016}, we utilized the CFRLM\_EN to model the transportation network. 
The CFRLM\_EN expands the original transportation network by adding a pseudo arc (an arc is a road segment between two adjacent transportation nodes in a path) between any two nodes in any path whose distance between each other is shorter than the PEVs' driving range. The scale of the binary charge decision variables of one path in the CFRLM\_EN is equal to the number of arcs (including the original arcs and the pseudo arcs)\footnote{Note that, to utilize the proposed service level model in this paper, the continuous charge decision variables in the CFRLM\_EN should be binary.}. By contrast, the scale of the binary charge choice variables of one path in the CFRLM\_SP is equal to the number of transportation nodes. Because that the number of transportation nodes is approximately equal to the number of original arcs, the sale of binary decision variables in the CFRLM\_SP is smaller than that in  the CFRLM\_EN. Interested readers can refer to \cite{Zhang_PlanFRLMTSG2016} for detailed introduction of the CFRLM\_EN. 

We compare the scales and solution time of Cases 1 and 3 (see Table~\ref{tab:case} for the parameter settings) utilizing the two different models in Table~{\ref{tab:scale}}. Note that the extra constraints introduced in Section \ref{sec:cut} are included in both of the two models. As expected, the scales of the binary variables of CFRLM\_SP are much smaller than those of the CFRLM\_EN. When assuming all the PEVs are homogeneous, the solution time of CFRLM\_EN is about fifteen times that of CFRLM\_SP. When the PEVs are categorized into four types, the CFRLM\_EN is intractable while the CFRLM\_SP can still be solved in less than 10 minutes.

\section{Conclusion}\label{section:conclusion}
{In this paper, we study the planning of PEV fast-charging stations on coupled transportation and power networks. We address three core questions in this problem: 1) how many charging spots should we construct in each station? 2) where should we locate these charging stations? 3) how do the transportation and power networks jointly impact the charging stations and the planning results? Specifically, first, we develop a closed-form service level model to describe a fast-charging station' service ability. This model can be used to determine the size of a  station servicing heterogeneous PEV charging demands. 
Then, we propose the modified CFRLM\_SP to explicitly capture time-varying PEV charging demands under driving range constraints in the transportation network. This model defines the feasible set of the PEV charging locations.
At last, we formulate a stochastic mixed-integer SOCP model to site and size fast-charging stations considering both the transportation and power network constraints.} 

{Numerical experiments validate the proposed methods. The simulation results show that the proposed service level model has high accuracy with heterogeneous PEV driving demands. The modified CFRLM\_SP can effectively describe the PEV driving range constraints considering time-varying traffic flow. It is also more computationally efficient than its counterpart in published literature. 
Simulation results also show that considering both the transportation constraints and the power network constraints with AC power flow at the same time leads to more economical investment decisions.}

In practice, the distribution networks may not be radial, as a result, the SOCP relaxation for the AC power flow model may not be exact. In that case, the planner may adopt alternative methods, e.g., the DC approximation used in our previous work \cite{Zhang_PlanFRLMTSG2016} or the semi-definite programming (SDP) relaxation \cite{OPF_Lavaei2012}. Adopting DC approximation for distribution systems in the planning will be computationally efficient. However, it may be less accurate (as discussed in Section \ref{sec:DC}). {Adopting SDP relaxation for power flow models in meshed networks may provide accurate solutions\footnote{{Note that the exactness of the SDP relaxation for AC power flow may not hold in some scenarios so that its solution may not be optimal or feasible for the original problem. See reference \cite{Low2014} for the detailed introduction.}}; however, it is generally less computationally efficient than the SOCP relaxation\footnote{{See reference \cite{Li2017} for an experimental analysis of the computational efficiency of both models.}}. The computational efficiency is predominant in the proposed planning model with a significant number of integer decision variables. However, to the best of our knowledge, there is no mature off-the-shelf commercial solver that can efficiently solve large-scale mixed-integer SDP problems.} 

Modeling the planning problem with meshed distribution networks are beyond the scope of this paper and will be our future work.
Computationally efficient solution method for the planning model in large-scale transportation and power network scenarios is also our future focus. 

\section*{Acknowledgment}
The authors would like to thank Prof. Zuojun Max Shen for insightful discussions.

\appendices
\section{Proof of Proposition 1}\label{proof_p1}
 In the following proof, we let index $e$ also denote the order of the arrival of a PEV, and PEV $e+y^\text{cs}$ is the $y^\text{cs}$st PEV that arrives after PEV $e$.

With [A2], $t_e^\text{c}=t_e^\text{a}$ is always satisfied, and we also have that $t_e^\text{d} = \min (t_e^\text{a}+T, t_{e+y^\text{cs}}^\text{a})$; therefore:
\begin{align}
&t_e^\text{d}-t_e^\text{a}\geq T \iff t_{e+y^\text{cs}}^\text{a}-t_e^\text{a}\geq T, \quad \forall e,\\
&\text{Pr}(t_e^\text{c}=t_e^\text{a}~\&~t_e^\text{d}-t_e^\text{c}\geq T) = \text{Pr}(t_{e+y^\text{cs}}^\text{a}-t_e^\text{a}\geq T), \forall e.
\end{align}

With [A1], the PEV arrival events are Poisson so that they are independent and identically distributed:
\begin{align}
&\text{Pr}(t_{e+y^\text{cs}}^\text{a}-t_e^\text{a}\geq T) = \text{Pr}(t_{y^\text{cs}}^\text{a}-t_0^\text{a}\geq T), \quad\forall e.
\end{align}

The probability that we observe the $y^\text{cs}$st PEV arrives after time $t_0^\text{a}+T$ is the same as the probability that we observe less than $y^\text{cs}$ PEVs that arrive from $t_0^\text{a}$ to $t_0^\text{a}+T$:
\begin{align}
\text{Pr}(t_{y^\text{cs}}^\text{a}-t_0^\text{a}\geq T)&= \text{Pr}(y^\text{ev}\leq y^\text{cs})\notag\\
&=\text{Pr}(t_e^\text{c}=t_e^\text{a}~\&~t_e^\text{d}-t_e^\text{c}\geq T),\quad \forall e.
\end{align}
This completes the proof. \hfill $\square$

\bibliographystyle{ieeetr}
\bibliography{ref}

\end{document}